# Pitchfork Bifurcation In A Coupled Cell System

*Shikhar Raj[1], Biplab Bose[1,*]*

[1]Department of Biosciences and Bioengineering, Indian Institute of Technology Guwahati, Guwahati, Assam, India 781039

*Email: biplabbose@iitg.ac.in

**Abstract:** Various biological phenomena, like cell differentiation and pattern formation in multicellular organisms, are explained using the bifurcation theory. Molecular network motifs like positive feedback and mutual repressor exhibit bifurcation and are responsible for the emergence of diverse cell types. Mathematical investigations of such problems usually focus on bifurcation in a molecular network in individual cells. However, in a multicellular organism, cells interact, and intercellular interactions affect individual cell dynamics. Therefore, the bifurcation in an ensemble of cells could differ from that for a single cell. This work considers a ring of identical cells. When independent, each cell exhibits supercritical pitchfork bifurcation. Using analytical and numerical tools, we investigate the bifurcation in this ensemble when cells interact through positive and negative coupling. We show that within a specific parameter zone, an ensemble of positively coupled cells behaves like a single cell with supercritical pitchfork bifurcation. In this regime, all cells are synchronized and have the same steady state. However, this unique behaviour is lost when cells interact through negative coupling. Apart from the synchronized (or homogenous) states, cell-cell coupling leads to certain heterogeneous steady states with unique patterns. We also investigate the distribution of such heterogeneous states under positive and negative coupling.

# 1. Introduction

Bifurcation theory is a mathematical framework that describes how small changes in one or more system parameters can lead to sudden qualitative changes in system behaviour. It is used to explain various biological phenomena, like cellular differentiation [1], cell cycle regulation [2] and epithelial-to-mesenchymal transition [3]. In these examples, molecular regulatory networks are considered as dynamical systems. These dynamical systems can have multiple steady states or fixed points. Bifurcation happens when changes in the values of one or more parameters change system behaviour by altering the number of steady states and/or their stability [4]. This shift marks a transition to a new dynamical regime and associated changes in cellular behaviour.

In developmental biology, bifurcation theory explains how identical cells diversify into different cell types. Diversification of cell types is a critical phenomenon in the embryonic development of multicellular organisms. Starting with a mass of identical cells, a whole organism made of different types of specialized cells emerges through repeated and sequential differentiation. Each cell type has unique characteristics (shape, size, location) and functions that define its phenotype.

Phenotypic diversification is regulated through molecular regulatory networks. These networks are dynamical systems with one or more stable steady states. In the dynamical systems paradigm of differentiation, these stable steady states correspond to different cell phenotypes [1]. During cellular differentiation, variations in gene expression or external cues make an existing stable steady state unstable and push the system towards a new stable steady state, leading to distinct cell types.

Such phenotypic diversifications have been explained by pitchfork and saddle-node bifurcations [1,5–7]. A specific cell type represents a particular stable branch in the bifurcation diagram. At the bifurcation point, new stable branches emerge, each representing a cell type. As the control parameter (environmental cue or expression of a regulator) crosses the bifurcation point, a cell moves to one of those stable branches.

Molecular networks regulating cellular processes are usually large and complicated. Even then, reduced dynamical models have been developed, and bifurcation theory has been used to explain cellular diversification [7–9]. Most of these works consider the dynamics of a network in a single, isolated cell. However, cells in a multicellular organism do not exist in isolation and have cell-cell communication. Cell-cell interaction can potentially alter the dynamics and bifurcation in a molecular network. In a multicellular organism, cells in an ensemble differentiate into different cell types, giving rise to unique cellular patterns. It is well established that cell-cell interactions play a crucial role in pattern formation during morphogenesis [10,11].

However, modelling the dynamics of a molecular network in an ensemble with cell-cell interactions is not trivial. As the number of cells increases, the dimension of the system grows enormously, limiting the scope of analytical and numerical investigations. Further, cell-cell interactions can be of many types, each with a distinct effect on the system dynamics and bifurcation.

These limitations are not unique to cellular systems. Generally, bifurcation analysis of an ensemble of dynamical systems with coupling/interaction is difficult. However, suppose individual dynamical systems are identical, and the coupled system has symmetry. In that case, one can use the concepts of equivariant dynamical systems to comment on the bifurcation behaviour of the whole system and the formation of patterns [12–14].

In this work, we consider a ring of identical cells. The molecular network in each cell has supercritical pitchfork bifurcation. Each cell interacts with its two neighbours, and cell-cell coupling affects the dynamics of the network. We considered two types of linear cell-cell coupling – positive and negative. As a whole, our systems of analysis have certain symmetries.

We investigate these coupled systems analytically and numerically. Our objective is to study steady state behaviours and check the supercritical pitchfork bifurcation in these coupled systems. We identify different types of synchronous (or homogenous) and nonsynchronous (or heterogenous) steady states that may arise and patterns that can form involving different cell types. By synchronous, we mean a steady state where all the cells have the same steady state. The bifurcation behaviour depends upon number of cells and the type of coupling. With many cells, the overall bifurcation diagram is very complicated. However, we show that with positive cell-cell coupling, the coupled system of any size has supercritical pitchfork bifurcation involving the synchronous steady states. More interestingly, a parameter zone exists where the system only has these synchronous steady states, and nonsynchronous ones do not exist. Therefore, in that parameter zone, this coupled system of $n$-cells behaves like a single-cell system with pitchfork bifurcation. However, this unique behaviour is lost with negative cell-cell coupling. We also show that negative cell-cell coupling is more effective than positive coupling in generating non-homogenous patterns.

## 2. The Model

We consider an array of *n* cells with periodic boundary conditions such that each cell has two neighbours. So, it is a ring of cells. Fig. 1a shows the arrangements for *n* = 3 and 6 cells.

We assume all cells have identical internal dynamics with supercritical pitchfork bifurcation. By internal dynamics, we mean the dynamics of a cell when it is isolated from others. To capture the internal dynamics of a cell, we use the normal form of supercritical pitchfork bifurcation [4],

$$\frac{dx_i}{dt} = \left(rx_i - x_i^3\right) \qquad x_i, r \in \mathbb{R} \tag{1}$$

Here, $x_i$ is the state variable for the *i*-th cell. We consider $x_i$ as a lumped variable that decides the phenotype of the *i*-th cell. The bifurcation parameter is *r*. This parameter is a proxy for the molecular signals that induce changes in phenotype.

To include the effect of interactions with the neighbours, Eq. (1) is modified,

$$\frac{dx_i}{dt} = \left(rx_i - x_i^3\right) + pG\left(x_{i-1}, x_{i+1}\right) \qquad x_i, r \in \mathbb{R} \tag{2}$$

$pG(x_{i-1}, x_{i+1})$ represents the neighbours' effect on the *i*-th cell dynamics. With periodic boundary conditions, $i+1=1$ when $i=n$ and $i-1=n$ when $i=1$.

Here, *p* is the coupling constant. We consider two models, one with positive coupling between cells ($p > 0$) and the other with negative coupling ($p < 0$).

$G(\cdot)$ is the coupling function. For simplicity, we considered an odd linear function that represents the average value of *x* in the neighbouring cells,

$$G(x_{i-1}, x_{i+1}) = \frac{x_{i-1} + x_{i+1}}{2}$$

So, the system of ODEs for an $n$-cell system is given by $n$ equations of the following form with the periodic boundary condition,

$$\frac{dx_i}{dt} = \left(rx_i - x_i^3\right) + p\left(\frac{x_{i-1} + x_{i+1}}{2}\right) \qquad x_i, r, p \in \mathbb{R}, \text{ and } i = 1, \cdots, n \tag{3}$$

### 3.1 Model behaviour without cell-cell interaction:

Without cell-cell coupling ($p = 0$), each cell undergoes independent supercritical pitchfork bifurcation (Fig. 1b). The bifurcation point is at $r = 0$. When $r < 0$, a cell has only one steady state $x_i^* = 0$, which is stable. The system, as a whole, achieves the stable steady state $x_1^* = x_2^* = \cdots = x_n^* = 0$. It is a homogenous or synchronous steady state. By homogeneous or synchronous, we mean an identical value of $x$ for all cells at the steady state.

Beyond the bifurcation point, the previous stable steady state turns unstable, and two new stable steady states ($x_i^* = \pm\sqrt{r}$) appear (Fig. 1b). As each cell is independent, an $n$-cell system has $2^n$ stable steady states. Out of these, only two are synchronous: $x_1^* = x_2^* = \cdots = x_n^* = \pm\sqrt{r}$.

### 3.2 Symmetry of the coupled system and its steady states

A dynamical system with symmetry can be investigated using concepts of group theory. Consider a dynamical system $\frac{d\mathbf{x}}{dt} = f(\mathbf{x})$, $\mathbf{x} \in \mathbb{R}^n$ with a group $G$ acting on its state space. The system is equivariant under $G$ if $f(g \cdot \mathbf{x}) = g \cdot f(\mathbf{x})$ for all $g \in G$ [14].

For a coupled system with $n$ cells, the system of ODEs in Eq. (3) has certain symmetries. When $n > 3$ this system is equivariant under cyclic group $\mathbb{Z}_n = \langle (x_1 x_2 \cdots x_{n-1} x_n) \rangle$. Here, $(x_1 x_2 \cdots x_{n-1} x_n)$ is a cyclic permutation of $\{x_1, x_2, \cdots, x_n\}$. When $n = 3$ this system is equivariant under the group $S_3$ of all permutations of $\{x_1, x_2, x_3\}$. In both cases, the system is also equivariant under $\mathbb{Z}_2 = \{I, -I\}$. Here, $I$ is an identity matrix.

Due to these symmetries, the coupled system has the following properties:

**a)** If $(x_1^*, x_2^*, \cdots, x_n^*)$ is a steady state, then $(-x_1^*, -x_2^*, \cdots, -x_n^*)$ is also a steady state, and they have the same stability.

**b)** For $n > 3$, if $(x_1^*, x_2^*, \cdots, x_n^*)$ is a steady state, then all cyclic permutations $(x_1^* x_2^* \cdots x_{n-1}^* x_n^*)$ are also steady states, and they have the same stability.

**c)** For $n = 3$, if $(x_1^*, x_2^*, x_3^*)$ is a steady state, then all permutations of $\{x_1^*, x_2^*, x_3^*\}$ are also steady states, and they have the same stability.

Note that these properties are valid for both positive and negative cell coupling.

*Proof of (a):* Let us represent our system of ODEs as $\dot{\mathbf{x}} = \mathbf{f}(x_1, x_2, \cdots, x_n)$.

At the steady state $(x_1^*, x_2^*, \cdots, x_n^*)$, $\mathbf{f}(x_1^*, x_2^*, \cdots, x_n^*) = 0$.

As the system is equivariant under $\mathbb{Z}_2 = \{I, -I\}$,

$$\mathbf{f}(-x_1^*, -x_2^*, \cdots, -x_n^*) = \mathbf{f}(-I(x_1^*, x_2^*, \cdots, x_n^*)) = -I \, \mathbf{f}(x_1^*, x_2^*, \cdots, x_n^*) = 0$$

Therefore, $(-x_1^*, -x_2^*, \cdots, -x_n^*)$ is also a steady state.

In this work, we analyze the local stability of fixed points by linearization. The Jacobian matrix of the system is,

$$\mathbf{J} = \begin{bmatrix} r - 3x_1^2 & \frac{p}{2} & 0 & \cdots & 0 & \frac{p}{2} \\ \frac{p}{2} & r - 3x_2^2 & \frac{p}{2} & \cdots & 0 & 0 \\ 0 & \frac{p}{2} & r - 3x_3^2 & \cdots & 0 & 0 \\ \vdots & \vdots & \vdots & \ddots & \vdots & \vdots \\ 0 & 0 & 0 & \cdots & r - 3x_{n-1}^2 & \frac{p}{2} \\ \frac{p}{2} & 0 & 0 & \cdots & \frac{p}{2} & r - 3x_n^2 \end{bmatrix}$$

**J** is a symmetric matrix. The off-diagonal elements are either $\frac{p}{2}$ or zero. For $n=3$, all off-diagonal elements are $\frac{p}{2}$. All $x_i$ appear only on the diagonal of the matrix and also in squared form. Therefore, the Jacobian matrix of $(x_1^*, x_2^*, \cdots, x_n^*)$ and $(-x_1^*, -x_2^*, \cdots, -x_n^*)$ are identical. So, the stability of both the steady states is the same. ∎

*Proof of (b)*: Our system is equivariant under $\mathbb{Z}_n = \langle (x_1 x_2 \cdots x_{n-1} x_n) \rangle$. Let $\sigma$ is the cyclic permutation $(x_1 x_2 \cdots x_{n-1} x_n)$. Let $(x_1^*, x_2^*, \cdots, x_n^*)$ be a steady state. The cyclic permutation of this steady state is $\sigma(x_1^*, x_2^*, \cdots, x_n^*) = (x_n^*, x_1^*, \cdots, x_{n-1}^*)$.

As our system is equivariant under $\sigma$,

$$\mathbf{f}(x_n^*, x_1^*, \cdots, x_{n-1}^*) = \mathbf{f}(\sigma(x_1^*, x_2^*, \cdots, x_n^*)) = \sigma \mathbf{f}(x_1^*, x_2^*, \cdots, x_n^*) = 0$$

Hence, $(x_n^*, x_1^*, \cdots, x_{n-1}^*)$ is also a steady state.

Let $\mathbf{J}^*$ be the Jacobian matrix for the steady state $(x_1^*, x_2^*, \cdots, x_n^*)$. The action of $\sigma$ could be achieved using a permutation matrix. Let that matrix be $\mathbf{P}$. Then, the Jacobian matrix of the steady state $(x_n^*, x_1^*, \cdots, x_{n-1}^*)$ would be

$$\mathbf{J}_p^* = \mathbf{P} \mathbf{J}^* \mathbf{P}^{-1} \tag{4}$$

$\mathbf{P}$ is an orthogonal matrix. Therefore, following Eq. (4), $\mathbf{J}^*$ and $\mathbf{J}_p^*$ are similar matrices and their eigenvalues must be identical. Therefore, the steady state $(x_1^*, x_2^*, \cdots, x_n^*)$ and its cyclic permutation $(x_n^*, x_1^*, \cdots, x_{n-1}^*)$ have the same stability. ∎

Proof of (c): The proposition for $n=3$ could be proved using logics similar to those used in the proof for proposition (b).

### 3.3 Synchronous steady states of the coupled system:

Let's assume this system has a homogenous or synchronous steady state $x_1^* = x_2^* = \cdots = x_n^* = \alpha$.

Using Eq. (3), $r\alpha - \alpha^3 + \dfrac{p}{2}(\alpha + \alpha) = 0$

$\therefore \alpha = 0$ or $\alpha = \pm\sqrt{(r+p)}$

When $r \leq -p$, the system has only one synchronous steady state at $x_1^* = x_2^* = \cdots x_n^* = 0$. When $r > -p$, there are three synchronous steady states, $x_1^* = x_2^* = \cdots = x_n^* = 0, \sqrt{r+p}, -\sqrt{r+p}$. Note that these synchronous steady states exist for both positive and negative coupling.

### 3.4 Bifurcation and synchronous steady states for positive coupling:

From the last section, we know that at $r = -p$ the number of synchronous steady states changes from one to three. It indicates bifurcation in the system. This section explores the bifurcation in the coupled system with positive coupling.

When $p > 0$, the coupled system has the following properties:

**a)** The system has a supercritical pitchfork bifurcation at $r = -p$ involving the synchronous steady states.

**b)** The nonsynchronous steady states are bounded by the synchronous ones.

*Proof for (a):* In the previous section, we have shown that for $r \leq -p$ the system has only one synchronous steady state $x_1^* = x_2^* = \cdots = x_n^* = 0$. However, for $r > -p$ it has two additional synchronous steady states $x_1^* = x_2^* = \cdots = x_n^* = \pm\sqrt{r+p}$.

Let, $r = -p + \delta$, where $\delta \to 0^+$. The steady states $x_i^* = \pm\sqrt{r+p} = \pm\sqrt{\delta} \to 0$. That means these two additional synchronous steady states originate from $x_1^* = x_2^* = \cdots = x_n^* = 0$ at $r = -p$. Therefore, $r = -p$ is a branch point.

We perform local stability analysis of these steady states through linearization. Consider a synchronous steady state $x_1^* = x_2^* = \cdots = x_n^* = \alpha$. The Jacobian matrix for this steady state will be a *n*-by-*n* circulant matrix,

$$\mathbf{J}^* = \begin{bmatrix} r-3\alpha^2 & \frac{p}{2} & 0 & \cdots & 0 & \frac{p}{2} \\ \frac{p}{2} & r-3\alpha^2 & \frac{p}{2} & \cdots & 0 & 0 \\ 0 & \frac{p}{2} & r-3\alpha^2 & \cdots & 0 & 0 \\ \vdots & \vdots & \vdots & \ddots & \vdots & \vdots \\ 0 & 0 & 0 & \cdots & r-3\alpha^2 & \frac{p}{2} \\ \frac{p}{2} & 0 & 0 & \cdots & \frac{p}{2} & r-3\alpha^2 \end{bmatrix}$$

For $n=3$, all off-diagonal elements are $\frac{p}{2}$.

As $\mathbf{J}^*$ is a circulant matrix, the *k*-th eigenvalue of this matrix is,

$$\lambda_k = \sum_{j=0}^{n-1} c_j \exp\left(\frac{2\pi jk}{n} i\right) \qquad (5)$$

Here $c_j$ are the elements of the first row of the matrix and $k = 0, 1, \cdots, n-1$.

For $\alpha = 0$,

$$\lambda_k = r + \frac{p}{2}\exp\left(\frac{2\pi k}{n} i\right) + \frac{p}{2}\exp\left(\frac{2\pi k(n-1)}{n} i\right) = r + p\cos\left(\frac{2\pi k}{n}\right) \qquad (6)$$

We know $-1 \le \cos\left(\frac{2\pi k}{n}\right) \le 1$. As $p > 0$, $p\cos\left(\frac{2\pi k}{n}\right) \le p$. Therefore, for $r < -p$, all eigenvalues are negative. When $r = -p$, one eigenvalue is zero ($\lambda_0 = 0$). The rest of the eigenvalues are negative. For $r > -p$ at least one eigenvalue is positive.

So, the steady state $x_1^* = x_2^* = \cdots = x_n^* = 0$ is stable for $r < -p$ and unstable for $r > -p$.

For $\alpha = \pm\sqrt{r+p}$, the *k*-th eigenvalue of $\mathbf{J}^*$ is

$$\lambda_k = -2r - 3p + p\cos\left(\frac{2\pi k}{n}\right) \tag{7}$$

As $p > 0$ and $r > -p$, all these eigenvalues are negative. Therefore, the synchronous steady states $x_1^* = x_2^* = \cdots = x_n^* = \pm\sqrt{r+p}$ are stable.

In essence, $r = -p$ is a bifurcation point. Before this point, this coupled system has only one synchronous and stable steady state. At the bifurcation point, one eigenvalue of the Jacobian matrix is zero. Beyond this point, the existing synchronous steady state changes its stability, and two new stable, synchronous steady states emerge. These changes in the number and stability of steady states are typical of supercritical pitchfork bifurcation [4,15]. Therefore, this coupled system of $n$ cells has supercritical pitchfork bifurcation involving the synchronous steady states. ∎

*Proof for (b):* Let $x_h > 0$ be the highest possible value of any $x$ when the steady state is nonsynchronous. Let $(x_1^*, \cdots, x_j^*, \cdots, x_n^*)$ is a nonsynchronous steady state with $x_j^* = x_h$.

Let, for this steady state, $x_j^* > x_{j+1}^*$ and $x_j^* > x_{j-1}^*$. So, $2x_j^* > x_{j-1}^* + x_{j+1}^*$. Using this inequality, we can write,

$$rx_j^* - (x_j^*)^3 + px_j^* > rx_j^* - (x_j^*)^3 + \frac{p}{2}(x_{j-1}^* + x_{j+1}^*)$$

As $(x_1^*, \cdots, x_j^*, \cdots x_n^*)$ is a steady state, using Eq. (3), $rx_j^* - (x_j^*)^3 + \frac{p}{2}(x_{j-1}^* + x_{j+1}^*) = 0$. Therefore, $rx_j^* - (x_j^*)^3 + px_j^* > 0$. Replacing $x_j^*$ by $x_h$ we get,

$$x_h\left[(r+p) - x_h^2\right] > 0 \tag{8}$$

As $x_h > 0$, then $(r+p) - x_h^2 > 0$. Therefore, when $p > 0$ and $r > -p$, $0 < x_h < \sqrt{r+p}$.

Following the system's symmetry, there must be another nonsynchronous steady state $(-x_1^*, \cdots, -x_j^*, \cdots, -x_n^*)$. As $x_h$ is the largest and $x_h > 0$, then $x_l = -x_h$ will be the smallest value in a nonsynchronous steady state such that $x_l > -\sqrt{r+p}$.

Therefore, when $p > 0$ and $r > -p$, the largest and smallest values of $x$ for a nonsynchronous steady state are bounded by the synchronous steady state values $\pm\sqrt{r+p}$.

The same bounds can be proved for any nonsynchronous steady state where $x_j^*$ is the largest or smallest.

Now consider $p > 0$ and $r \leq -p$. In this case, the inequality $(r+p) - x_h^2 > 0$ in Eq. (8) does not hold. That means for $p > 0$ and $r \leq -p$, the system cannot have a nonsynchronous steady state with $x_h > 0$. Using the logic of symmetry, the system also cannot have a lowest value $x_l < 0$. Therefore, when $r \leq -p$ the system cannot have a nonsynchronous steady state.

Hence, for all values of $r$, the nonsynchronous steady states are bounded by the synchronous steady states. ∎

The number of nonsynchronous steady states depends on the number of cells, and it is difficult to find them analytically for an arbitrary value of $n$. Therefore, we numerically investigated the steady states and bifurcations for several values of $n$.

We have shown that at the bifurcation point $r = -p$, one eigenvalue of the Jacobian matrix for the synchronous steady state $x_1^* = x_2^* = \cdots = x_n^* = 0$ turns zero, and two new synchronous steady states emerge. Following Eq. (6), two eigenvalues of this Jacobian matrix, $\lambda_1$ and $\lambda_{n-1}$, become zero at $r = -p\cos\left(\frac{2\pi}{n}\right)$. Through numerical analysis, we observed that multiple nonsynchronous steady states emerge at this position. Nonsynchronous steady states do not exist when $r \leq -p\cos\left(\frac{2\pi}{n}\right)$.

Therefore, when $r \leq -p\cos\left(\frac{2\pi}{n}\right)$, this system with $n$ cells has only synchronous steady states. The pitchfork bifurcation involves those synchronous steady states. Therefore, when $r \leq -p\cos\left(\frac{2\pi}{n}\right)$, the system of $n$-cells with positive coupling behaves like a one-dimensional system with supercritical pitchfork bifurcation. In this regime, the steady state behaviour of the $n$-cell system is represented by the ODE,

$$\frac{dx}{dt} = (r+p)x - x^3$$

Numerical investigation showed that, the nonsynchronous steady states emerging from $r = -p\cos\left(\frac{2\pi}{n}\right)$ are all unstable. Numerical analysis showed that stable nonsynchronous steady states appear further away from $r = -p\cos\left(\frac{2\pi}{n}\right)$. Therefore, there is an extended zone beyond $r = -p\cos\left(\frac{2\pi}{n}\right)$, where all stable steady states are synchronous, and the system effectively behaves like a one-dimensional system.

This behaviour can be easily understood from the bifurcation diagram for a three-cell system shown in Fig. 1c. We used MatCont [16] for numerical bifurcation analysis. Fig. 1c shows the bifurcation diagram for $n = 3$ and $p = 0.5$. The diagram shows the steady state behaviour of $x_1$ for different $r$ values. The bifurcation diagrams for $x_2$ and $x_3$ are identical.

As expected, the system has a supercritical pitchfork bifurcation at $r = -p = -0.5$. Before this bifurcation point, the system has only one steady state $(0,0,0)$. It is a synchronous steady state and stable. At the bifurcation point, one eigenvalue of the Jacobian matrix for $(0,0,0)$ turns zero. Beyond this bifurcation point, the steady state $(0,0,0)$ is unstable. At the bifurcation point, two stable synchronous steady states $\left(\sqrt{r+p}, \sqrt{r+p}, \sqrt{r+p}\right)$ and $\left(-\sqrt{r+p}, -\sqrt{r+p}, -\sqrt{r+p}\right)$ appear.

As $r$ increases, the steady state $(0,0,0)$ undergoes further changes in its behaviour. At $r = -p\cos\left(\frac{2\pi}{n}\right) = 0.25$ two eigenvalues of the Jacobian matrix for (0, 0, 0) are zero, and one is positive. Beyond this point, all the eigenvalues are positive, and the steady state $(0,0,0)$ is a repeller.

Numerical analysis shows that twelve nonsynchronous, unstable steady states emerge at $r = 0.25$. Six of these are of the form $(0, \alpha, -\alpha)$ and rest are of the form $(\alpha, \beta, \alpha)$. Note that all these steady states are not visible in the bifurcation diagram (Fig. 1c) as those are overlapping.

In Fig. 1c, six saddle-node bifurcations appear at $r \approx 1.35$. Three are above the $(0,0,0)$ branch, and the rest are below it. Note that two such bifurcations are not visible in Fig. 1c as those are overlapping. The steady states involved in these bifurcations are of the form $(\alpha, \beta, \alpha)$. As derived earlier, all these nonsynchronous steady states are bounded by the synchronous steady states (Fig. 1c).

In this work, we are primarily concerned about the stable steady states as they represent particular phenotypes. The bifurcation diagram in Fig. 1c has three distinct zones of stable steady states. For $r \leq -0.5$ the system has only one stable steady state that is synchronous. For $-0.5 < r \leq 1.35$ the system has two stable steady states that are synchronous. Therefore, up to $r = 1.35$ this three-cell system effectively behaves like a one-cell system with supercritical pitchfork bifurcation. Beyond this point, the system has eight stable steady states involving synchronous and nonsynchronous ones.

We performed bifurcation analysis with different values of $r$ and $p$. The bifurcation pattern in Fig. 1c was observed for all those $r$-$p$ combinations. Fig. 1d shows the number of stable steady states for different combinations of $r$ and $p$.

Fig. 2a and 2b show the number of stable steady states for different combinations of $r$ and $p$, for $n = 4$ and 6. Notice that all three $r$-$p$ plots (Fig. 1d and Fig. 2) have distinct zones for one and two stable steady states. These steady states are synchronous. Within these two zones, the $n$-dimensional ensemble behaves like a one-dimensional system with supercritical pitchfork bifurcation. For higher values of $r$, the systems have zones with more stable steady states. In these zones, only two stable steady states are synchronous; the rest are nonsynchronous.

### 3.5 Bifurcation and synchronous steady states for negative coupling:

Symmetries in the system and steady states are valid for both positive and negative coupling. Similarly, the same synchronous steady states exist for positive and negative coupling (refer to section 3.3). $r = -p$ is a branch point even when $p < 0$. For $r < -p$ the system has only one synchronous steady state $x_1^* = x_2^* = \cdots = x_n^* = 0$. When $r > -p$ it has two additional synchronous steady states $x_1^* = x_2^* = \cdots = x_n^* = \pm\sqrt{r+p}$.

However, the stability of these steady states and the bifurcations in the system are distinctly different when $p < 0$. For negative coupling, the system has the following properties:

**a)** The nonsynchronous steady states are not bounded by the synchronous steady states.

**b)** The bifurcation involving the synchronous steady states is not supercritical pitchfork bifurcation.

*Proof for (a):* Let $x_h > 0$ be the highest possible value of any $x$ when the steady state is nonsynchronous. Let $(x_1^*, \cdots, x_j^*, \cdots, x_n^*)$ be a nonsynchronous steady state with $x_j^* = x_h$, $x_j^* > x_{j+1}^*$ and $x_j^* > x_{j-1}^*$. So, $2x_j^* > x_{j-1}^* + x_{j+1}^*$. Using this inequality, we can write,

$$rx_j^* - (x_j^*)^3 + px_j^* < rx_j^* - (x_j^*)^3 + \frac{p}{2}(x_{j-1}^* + x_{j+1}^*) \quad \text{(note that in this case } p < 0\text{)}$$

As $(x_1^*, \cdots, x_j^*, \cdots x_n^*)$ is a steady state, using Eq. (3), $rx_i^* - (x_i^*)^3 + px_i^* < 0$. Replacing $x_j^*$ by $x_h$ we get,

$$x_h\left[(r+p) - x_h^2\right] < 0 \tag{9}$$

As $x_h > 0$, then $(r+p) - x_h^2 < 0$. Therefore, when $p < 0$ and $r > -p$, $x_h > \sqrt{r+p}$.

Following the system's symmetry, there must be another nonsynchronous steady state $(-x_1^*, \cdots, -x_j^*, \cdots, -x_n^*)$. If $x_h$ is the largest value of $x$, then $x_l = -x_h$ will be the smallest value of $x$ in a nonsynchronous steady state such that $x_l < -\sqrt{r+p}$.

Hence, unlike the system with positive coupling, the nonsynchronous steady states are not bounded by the synchronous steady states. Rather, for $p < 0$, the synchronous steady states are bound by one or more of the nonsynchronous steady states. ∎

*Proof for (b):* The Jacobian matrices for the synchronous steady states and the formulations for their eigenvalues remain the same, as discussed in section 3.4. Let us first analyze the local stability of the synchronous steady state $x_1^* = x_2^* = \cdots = x_n^* = \pm\sqrt{r+p}$.

The $k$-th eigenvalue of the Jacobian matrix $\mathbf{J}^*$ is

$$\lambda_k = -2r - 3p + p\cos\left(\frac{2\pi k}{n}\right)$$

As $p < 0$, the largest eigenvalue is obtained for the smallest $\cos\left(\frac{2\pi k}{n}\right)$. Therefore, when $n$ is even, and $n \geq 3$ the largest eigenvalue is,

$$\lambda_{max} = -2r - 3p - p = -2(r + 2p)$$

The largest eigenvalue is negative when $r > -2p$. Therefore, for $r > -2p$, all eigenvalues of the Jacobin matrix are negative, and these synchronous steady states are stable. Otherwise, they are unstable.

When, $n$ is odd and $n \geq 3$, the minimum value of $\cos\left(\frac{2\pi k}{n}\right)$ is obtained for $k = \frac{n+1}{2}$ and $k = \frac{n-1}{2}$.

So, the largest eigenvalue is,

$$\lambda_{max} = -2r - 3p + p\cos\left(\frac{\pi(n \pm 1)}{n}\right)$$

To make all the eigenvalues negative, $\lambda_{max}$ must be negative. However, that depends upon $n$. When $n \gg 1$, $\lambda_{max} \approx -2r - 3p - p$. In such a case, the synchronous steady states $x_1^* = x_2^* = \cdots = x_n^* = \pm\sqrt{r+p}$ will change stability from unstable to stable near $r = -2p$. For $n = 3$, this change in stability happens at $r = -1.75p$.

Now, we analyze the stability of the other synchronous steady state $x_1^* = x_2^* = \cdots = x_n^* = 0$. In this case, the $k$-th eigenvalue of the Jacobian matrix $\mathbf{J}^*$ is

$$\lambda_k = r + p\cos\left(\frac{2\pi k}{n}\right)$$

The maximum eigenvalue is $\lambda_{max} = r - p$ when $n$ is even (for $k = \frac{n}{2}$). All the eigenvalues will be negative if $\lambda_{max} = r - p < 0$. Therefore, when $r < p$ this synchronous steady state is stable. For $r > p$, at least one eigenvalue will be positive and this steady state will be unstable.

When $n$ is odd, $\lambda_{max}$ will depend upon $n$. It can be shown that when $n$ is odd and $n \gg 1$, $\lambda_{max} \approx r - p$. Therefore, the stable to unstable transition will happen near $r = p$. For $n = 3$ this transition happens at $r = 0.5p$.

In essence, the stabilities of synchronous steady states do not match with the supercritical pitchfork bifurcation. So, for the system with negative coupling, the bifurcation involving the synchronous steady states is not supercritical pitchfork bifurcation. ∎

Fig. 3a shows the bifurcation diagram for a three-cell system with negative coupling ($p = -0.5$). The diagram shows the steady state behaviour of $x_1$ for different $r$ values. The bifurcation diagrams for $x_2$ and $x_3$ are identical. The branches for synchronous steady states are marked by S. These branches have the shape of pitchfork bifurcation, with the branch point at $r = -p = 0.5$. However, this bifurcation is not a supercritical pitchfork. The steady state $x_1^* = x_2^* = x_3^* = 0$ changes from stable to unstable at $r = 0.5p$, not at $r = -p$. Further, the branches for $x_1^* = x_2^* = x_3^* = \pm\sqrt{r+p}$ are initially unstable and turn stable near at $r = -1.75p$.

At the branch point $r = -0.5p$, twelve nonsynchronous steady states emerge (Fig. 3a). Note that all these steady states are not visible in the bifurcation diagram as those are overlapping. Out of these, six are stable and are of the form $(\alpha, \beta, \alpha)$. The positively coupled system also has a branch point at this position, but all the steady states emerging from that point are unstable.

The bifurcation diagram with negative coupling has two other differences with the one with positive coupling. Unlike Fig. 1c, the branches of synchronous steady states in Fig. 3a are bounded by the nonsynchronous one. We have already proved the same analytically earlier in this section.

Like the positively coupled system, this system also has limit points. However, the steady states involved in these limit points are unstable. Therefore, we do not have saddle-node bifurcation that is observed for positive *p* (Fig. 1c).

The number of stable steady states for various combinations of *r* and *p* for three-cell and four-cell systems are shown in Fig. 3b and 3c, respectively. The zone of one stable steady state in both plots represents the synchronous steady state $x_1^* = x_2^* = \cdots = x_n^* = 0$.

For *n* = 3, there is no zone of two stable steady states. The stable synchronous steady states $x_1^* = x_2^* = x_3^* = \pm\sqrt{r+p}$ appear in the zone with eight stable steady states.

For *n* = 4, there is a zone of two stable steady states, but those steady states are nonsynchronous. They have the form of $(-\alpha, \alpha, -\alpha, \alpha)$. The stable synchronous steady states appear in the zone with six stable steady states.

In essence, with negative coupling, the system does not have a parameter zone where it behaves like a one-dimensional system with supercritical pitchfork bifurcation. It also does not have a parameter zone with only non-zero, synchronous, stable steady states.

### 3.6. Diversification of cell types and pattern formation

Till now, we primarily investigated the synchronous steady states. In this section, we focus on the nonsynchronous steady state and the formation of different patterns of cell states. Each stable steady state represents one cell type/state, and bifurcation in the system leads to the diversification of cell types. When the system is at a synchronous steady state, all cells are in the same state, and the pattern will be homogenous. However, non-homogenous patterns arise when the system is in a nonsynchronous steady state.

The number of nonsynchronous steady states depends upon the number of cells in the system. We numerically investigated the pattern formation with a specific number of cells. For brevity, we discuss the results of a 4-cell system only.

We used the DifferentialEquations.jl package [17] of Julia [18] for simulations. The steady state depends upon the initial values and values of the parameters *r* and *p*. We performed simulations for various combinations of *r* and *p*. For each case, 10,000 independent simulations were performed with randomly chosen initial values spanning a large state space covering all possible steady states. Subsequently, the steady state behaviours of the system were evaluated, and distributions of different cell state patterns were calculated.

Results for a four-cell system are shown in Figure 4. In this case, the coupling is either positive ($p = 1$) or negative ($p = -1$), and the $r$ varies. For a low $r$-value ($r = 0.2$), the positively coupled system is in the zone with only synchronous steady states. Therefore, the system reaches two homogenous patterns with equal frequency in our numerical simulations (Fig. 4a). However, for this value of $r$, a negatively coupled system has only nonsynchronous steady states with the pattern ($a, -a, a, -a$).

As $r$ is increased to one, a nonsynchronous pattern appears ($-a, -a, a, a$) for both the positive and negative coupling. Even then, the positively coupled system reaches a synchronous pattern for ~80% of initial values (Fig. 4b). For the negatively coupled system ($a, -a, a, -a$) is the dominant pattern (Fig. 4b).

For the positively-coupled system, the synchronous steady states are also dominant for $r = 2.5$ (Fig. 4c). However, nonsynchronous steady states dominate when $r$ is very high ($r = 4$ in Fig. 4d). In fact, some nonsynchronous patterns appear with the same frequency for both positive and negative coupling. For both $r = 2.5$ and 4, the negatively coupled system has homogenous patterns, but those appear in a minority of cases.

These simulations show that the non-homogenous patterns are favoured in the negatively coupled system. However, those patterns dominate only for high $r$ in the positively coupled system.

**4. Coupled mutual repressors:**

In the model defined in section 2, we have used the normal form of supercritical pitchfork bifurcation to represent the internal dynamics of individual cells. This section considers a mutual repressor circuit in each cell that shows supercritical pitchfork bifurcation. Molecular circuits with mutual repression are involved in cell differentiation [19,20].

The following system of ODEs is used for the mutual repressor circuit in each cell,

$$\frac{dx_i}{dt} = \frac{r}{1+y_i^2} - x_i$$
$$\frac{dy_i}{dt} = \frac{r}{1+x_i^2} - y_i \quad ; \quad x_i \geq 0, y_i \geq 0, r \geq 0 \tag{10}$$

Here, we can consider *x* and *y* as two genes crucial in a particular phenotypic diversification or cell fate-determining molecules. So, $(x_i, y_i)$ represents the state of a cell.

This system of ODEs has supercritical pitchfork bifurcation, with *r* being the bifurcation parameter (Fig. 5a). Before the bifurcation point, the system has a stable steady state $x_i^* = y_i^*$. Two new stable branches emerge at the bifurcation point, and the existing $x_i^* = y_i^*$ branch turns unstable. However, unlike the normal form of pitchfork bifurcation, these two stable branches are not symmetric around the central branch. In general, this system of ODEs (10) does not have the $x \to -x$ symmetry of the normal form of pitchfork bifurcation.

Fig. 5a shows the bifurcation diagram of $x_i$. The bifurcation diagram of $y_i$ is the same but has a crucial difference. In the bistable zone, $x_i^* \neq y_i^*$. If the system is on the upper branch of the bifurcation diagram of *x*, it will be in the lower branch in the bifurcation diagram of *y*. That means a cell with high *x* will have low *y* and *vice versa*.

Like the model in section 2, we consider a ring of *n*-cells. So, the system of ODEs for an *n*-cell system is given by 2*n* equations of the following forms,

$$\frac{dx_i}{dt} = \frac{r}{1+y_i^2} - x_i + p\left(\frac{x_{i-1} + x_{i+1}}{2}\right)$$
$$\frac{dy_i}{dt} = \frac{r}{1+x_i^2} - y_i + p\left(\frac{y_{i-1} + y_{i+1}}{2}\right) \tag{11}$$

Like in Eq. (3), *p* in Eq. (11) is the coupling constant representing either positive or negative coupling.

For $n > 3$, this system is equivariant under the cyclic group $\mathbb{Z}_n = \langle (s_1 s_2 \cdots s_{n-1} s_n) \rangle$, where $s_i = (x_i, y_i)$ for *i* = 1 to *n*. When $n = 3$, this system is equivariant under the group $S_3$ of all permutations of $\{s_1, s_2, s_3\}$. However, this system does not have the reflection symmetry observed in our previous model.

We investigated this model numerically for different numbers of cells. As we have two state variables for each cell, we define a synchronous steady state as $\left(x_1^* = x_2^* = \cdots = x_n^*, y_1^* = y_2^* = \cdots = y_n^*\right)$. There can be two types of synchronous steady states − either $x_i^* = y_i^*$ or $x_i^* \neq y_i^*$ for all $i$.

Fig. 5b shows the bifurcation diagram for a three-cell system with a positive coupling. This system has supercritical pitchfork bifurcation involving synchronous steady states. There are no nonsynchronous steady states.

Before the bifurcation, there is only one stable steady state that is synchronous− $\left(x_1^* = x_2^* = x_3^* = y_1^* = y_2^* = y_3^*\right)$. In the bistable zone, the steady states are also synchronous. However, in two synchronous, stable steady states $x_i^* \neq y_i^*$. As explained earlier, in this zone, the behaviour of $x_i$ is the opposite of $y_i$. When $x_i$ is on the upper stable steady state branch, $y_i$ is on the lower stable steady state branch.

We performed the bifurcation analysis for various combinations of $r$ and $p$ (Fig. 5d). The parameter space, $r \geq 0$ vs. $0 \leq p < 1$, has two zones. In one zone, there is only one synchronous stable steady state of the form $x_i^* = y_i^*$. The other zone has two synchronous, stable steady states of the form $x_i^* \neq y_i^*$. Numerical analysis shows that for $p > 1$, steady states values are negative and unstable.

So, with suitable positive cell-cell coupling ($0 \leq p < 1$), this system behaves like a single-cell system with supercritical pitchfork bifurcation. In this regime, the steady state behaviour system is by a system of ODEs of two variables,

$$\frac{dx}{dt} = \frac{r}{1+y^2} + (p-1)x$$
$$\frac{dx}{dt} = \frac{r}{1+x^2} + (p-1)y$$

However, a negatively coupled system does not have this behaviour. Fig. 5c shows the bifurcation diagram for the three-cell system with $p = -0.5$. Now, the system has both synchronous and nonsynchronous steady states, and synchronous steady states are bounded by nonsynchronous ones. Note that some of the steady states are negative; hence, they are not biologically plausible.

The synchronous branches in this diagram have the shape of pitchfork bifurcation, with the branch point at $r = 3$. However, this bifurcation is not a supercritical pitchfork. The branch for $x_i^* = y_1^*$ changes from stable to unstable at $r = 1$ not at 3. $r = 1$ is also a branch point. Several new nonsynchronous steady states emerge from this point. Further, the two new branches emerging at $r = 3$ are initially unstable and subsequently turn to stable near $r = 6$.

Fig. 5e shows the number of stable steady states for different combinations of $r$ and $p$ when cells are negatively coupled. There is no zone of two stable steady states. Cells are synchronized in the zone of one stable steady state. There is no synchronization in the zone of six stable steady states. Two synchronized states appear in the zone with eight stable steady states.

In essence, the bifurcation behaviour of this system is largely similar to the model using the normal form of supercritical pitchfork bifurcation. This similarity was observed even in systems with a higher number of cells (data not shown to keep the manuscript concise).

## 5. Discussions:

In the present work, we investigated the bifurcation in coupled cellular systems, where each cell has internal dynamics with supercritical pitchfork bifurcation. Cells are coupled through positive and negative cell-cell interactions, and we investigated the bifurcation of these multicellular systems. Our interest lies in finding synchronization where a multicellular system behaves like a single cell.

We show that positive and negative coupling between cells have distinct effects. When cells are positively coupled, a region exists in the parameter space where the multicellular system is effectively one cell and shows supercritical pitchfork bifurcation. In that parameter space, cells are always synchronized. This unique behaviour is lost when the coupling is negative.

We investigated this behaviour using two models. In the first model, the internal dynamics of each cell is represented by the normal form of the supercritical pitchfork bifurcation. The state variables do not represent any particular molecule but can be considered as lumped variables that decide the cell fate. The advantage of this model is that it can be investigated using analytical methods. Subsequently, we used another model where each cell has a mutual repressor circuit involving two genes. As analytical investigations are difficult, we investigated this model numerically. Interestingly, the key observations on cell synchronization and supercritical pitchfork bifurcation are the same in both models.

In biological systems, adjacent cells get coupled through cell surface molecules. Such lateral interactions are crucial in the differentiation of cells and pattern formation in multicellular systems. The Delta-Notch signalling is a well-known example of lateral interactions between cells [21,22]. The Delta ligand on one cell binds to the Notch receptor on its adjacent cell, activating the Notch pathway in the neighbour. This leads to the suppression of Delta expression in this neighbour cell. Such lateral inhibition leads to a salt-and-pepper pattern, with one cell having a high Notch signalling (receiver phenotype), whereas its neighbour has a low Notch activity (sender phenotype).

Lateral induction is also possible in Notch signalling. This involves Jagged, another Notch ligand. Jagged on a cell binds to the Notch receptor on its adjacent cell, activating Notch signalling. This leads to upregulation in the expression of Jagged in this neighbour cell. Eventually, such lateral induction leads to a homogenous state with all cells having the same phenotype [21,23].

Lateral inhibition in the Notch signalling is equivalent to negative cell-cell coupling in our models. Lateral inhibition in Delta-Notch signalling leads to the formation of patterns with cells in different states. Similarly, negative coupling in our model systems favours nonsynchronous steady states and the formation of diverse patterns.

On the other hand, lateral induction in the Notch signalling is equivalent to positive coupling in our models. Like Jagged-Notch lateral induction, positive coupling in our model systems favours synchronization with the whole system attaining a homogenous state. However, heterogeneous patterns dominate when the bifurcation parameter $r$ is high.

Here, we use the examples from Delta-Notch signalling to highlight the existence of different types of cell-cell interactions in multicellular systems. That does not mean that our models represent Delta-Notch signalling. However, equations similar to the normal form of pitchfork bifurcation have been used to model pattern formation in the Delta-Notch system [24,25], and pitchfork bifurcations have been reported in several models of the Delta-Notch system [26–28].

In our models, we considered a system of multiple interacting cells. However, our questions and observations are equally relevant for any arbitrary dynamical system with interacting sub-systems. Galizia and Piiroinen [29] investigated a similar question for a network of connected identical agents. They wanted to find the condition under which the dynamical properties of the entire network are equivalent to those of the individual agents. They defined the Region of Reduced Dynamics (RRD) as the subset of the parameter space where the steady states of uncoupled agents and the stability of those steady states are retained in the coupled network. Therefore, in the RRD, the system with $n$-agents is effectively a system of one agent.

This is similar to what we have observed for a positively coupled system in our work. However, the synchronization in our analysis is less restrictive. In the RRD, agents are not only synchronized, but steady states in the coupled and uncoupled systems are the same. In our work, the bifurcation of the coupled system is similar to that of an individual cell. However, the steady states are not identical in coupled and uncoupled cases.

In our models, when cells are coupled positively, a region of the parameter space exists where the steady state behaviour of the coupled system is similar (but not identical) to that of an uncoupled cell. Therefore, this region could also be regarded as a region of reduced dynamics in a broader sense.

For our first model, this region is just a small subset of the whole parameter space. However, for the second model, with positive coupling, the whole parameter space allows synchronization, and the system dynamics effectively reduces to one-cell dynamics.

We analyzed only two types of system of ODEs that show supercritical pitchfork bifurcation. Many other systems of ODEs can show supercritical pitchfork bifurcation. Our present work has not investigated whether our observations are valid for a ring of agents with any system of ODEs having supercritical pitchfork bifurcation. Also, changes in the number of neighbours and the rules of interactions may invalidate observations made in this work. It will be interesting to expand this work further to make more generalized statements about supercritical pitchfork bifurcation in coupled systems.

Figure 1

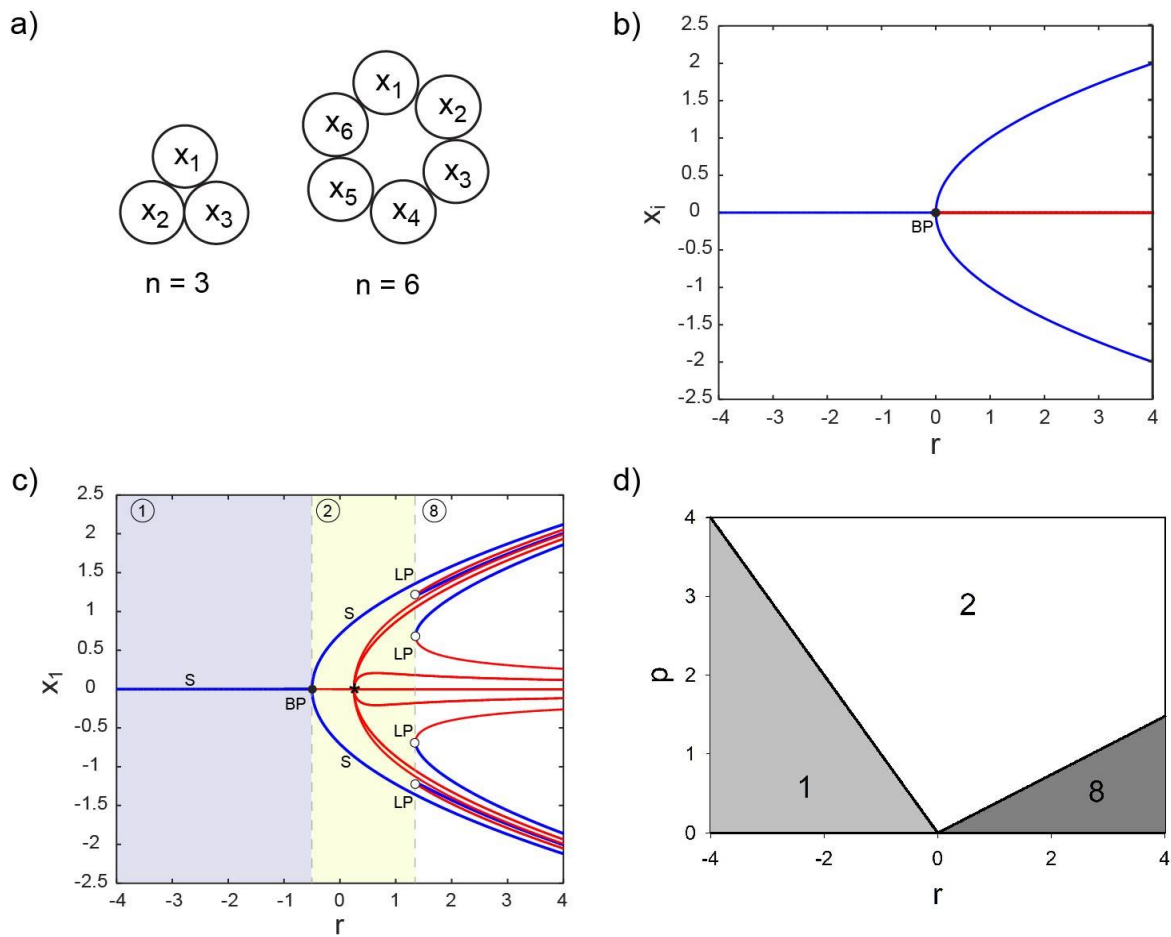

Figure 1: Bifurcation in a ring of cells. a) shows the arrangement for three and six cells. Each cell has two neighbours. b) The bifurcation diagram for $x$ when cells are not coupled ($p = 0$). c) The bifurcation diagram for a three-cell system with positive cell-cell coupling ($p = 0.5$). In (b) and (c), stable and unstable branches are blue and red, respectively. In (c), the branches for synchronous stable steady states are marked S. This diagram is divided into multiple zones based on the number of stable steady states of the whole system (written within circles). BP: Branch Point. LP: Limit Point. d) Phase diagram of the three-cell system showing the number of stable steady states for different $r$ and $p$ values. It has three zones that are coloured differently. The number of stable steady states in each zone is written on the diagram.

Figure 2:

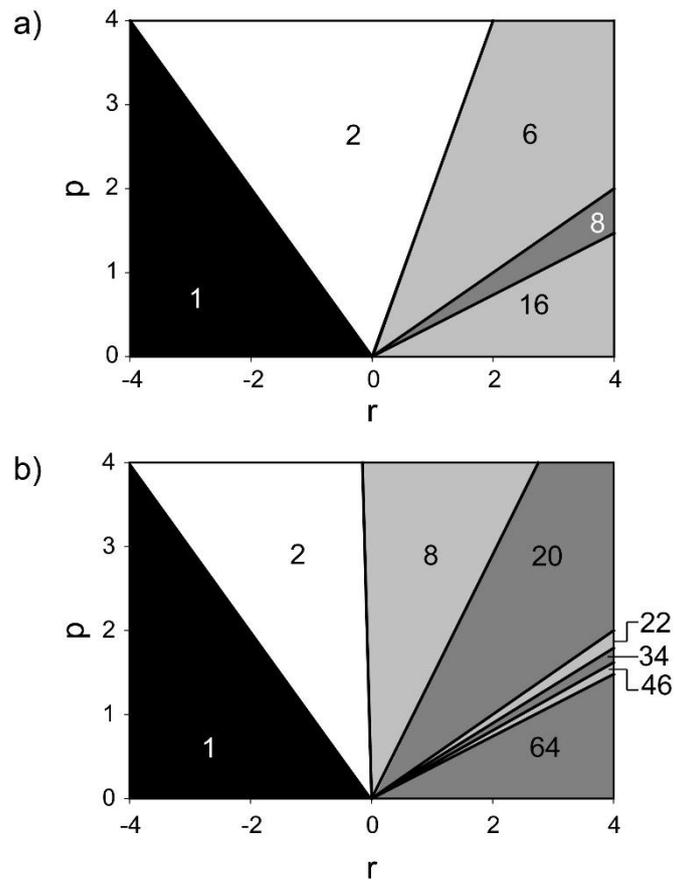

Figure 2: Phase diagram of a) four-cell and b) six-cell systems with positive coupling showing the effect of *r* and *p* on the number of stable steady states. Each diagram is divided into zones that are shaded with different colours. The number of stable steady states in each zone is written on the diagram.

Figure 3.

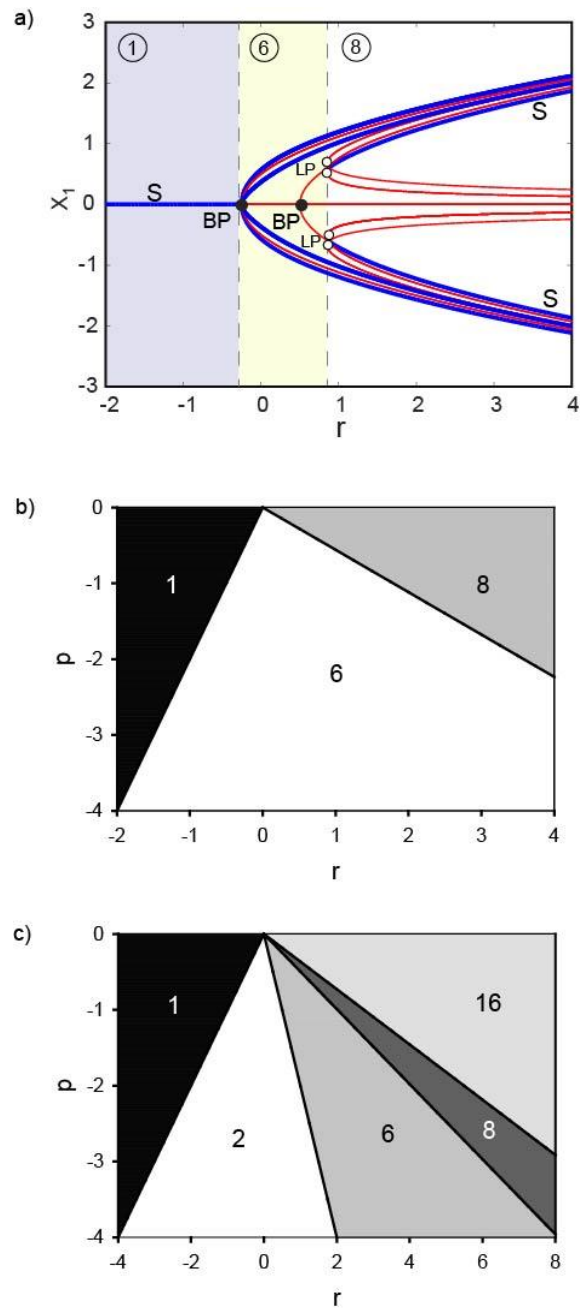

Figure 3. Behavior of the system with negative cell-cell coupling. a) The bifurcation diagram for a three-cell system with negative cell-cell coupling ($p = -0.5$). Stable and unstable branches are blue and red, respectively. The branches for synchronous, stable steady states are marked by S. BP: Branch Point. LP: Limit Point. This diagram is divided into multiple zones based on the number of stable steady states of the whole system (written within circles). b) and c) are the phase diagrams for $n = 3$ and $n = 4$, showing the number of stable steady states for different $r$ and $p$ values. Zones are marked with different colours. The number of stable steady states in each zone is written on the diagram.

Figure 4.

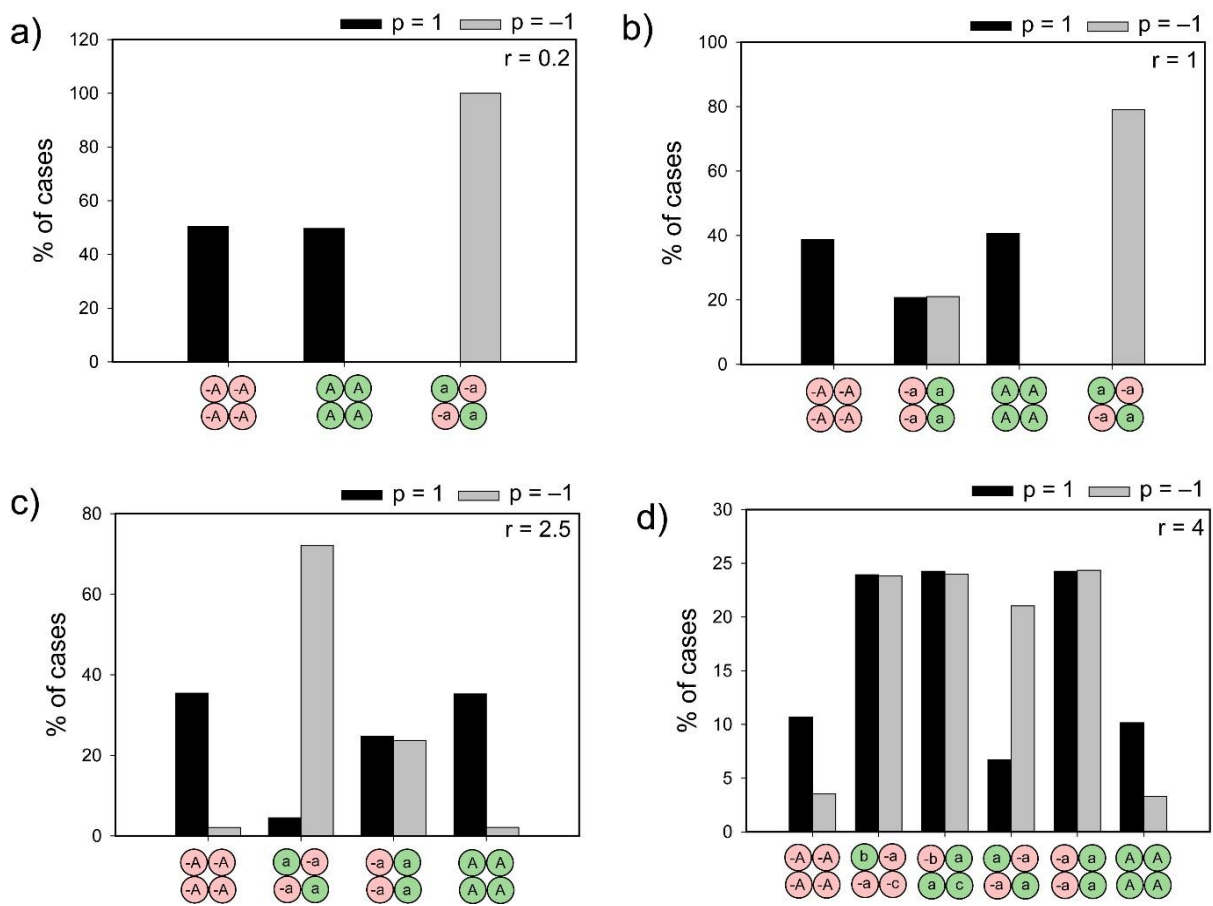

Figure 4: Patterns of homogenous and heterogenous cell states and their distributions in a four-cell system. The horizontal axes represent different steady state patterns. The vertical axis represents the percentage of simulations that reached a particular pattern. Cells are represented by circles. In all cases, red represents a negative steady state value of *x*, and the green represents a positive value. Steady state values are represented symbolically within the circles. A and −A represent the synchronous steady state values $\sqrt{r+p}$ and $-\sqrt{r+p}$, respectively. The lowercase alphabets represent other steady state values. However, they do not represent any specific numerical value and are used only to represent patterns.

Figure 5.

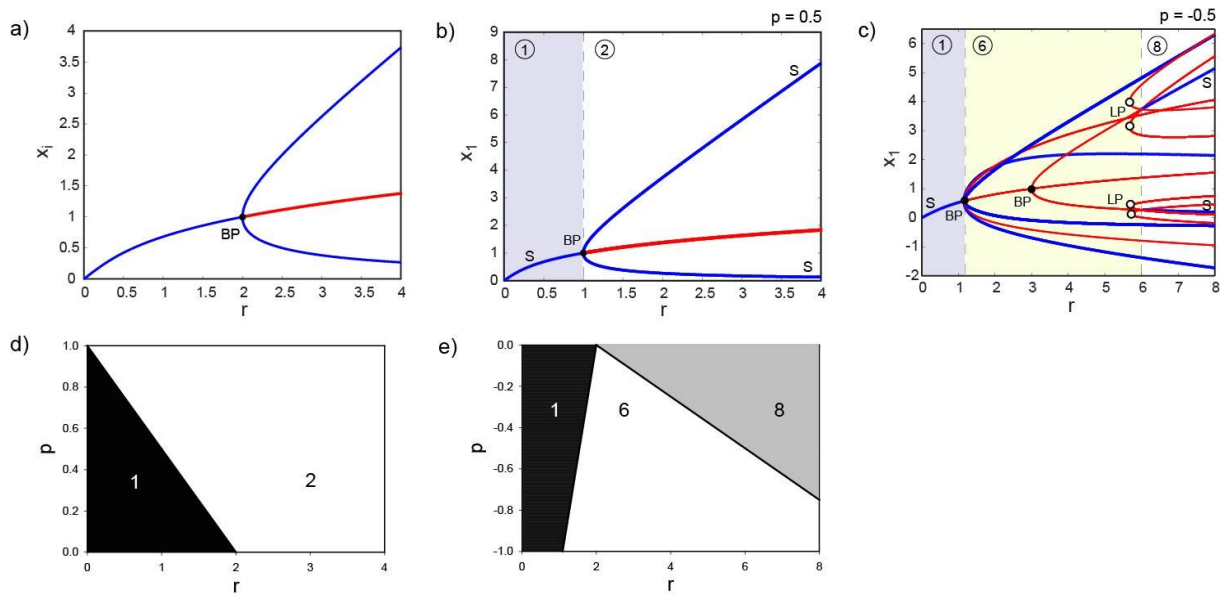

Figure 5: Bifurcations in a ring of cells with mutual repressor. a) The bifurcation diagram for $x$ in a cell when cells are uncoupled. b) and (c) show bifurcation diagrams for a three-cell system with positive ($p = 0.5$) and negative cell-cell coupling ($p = -0.5$), respectively. In all the bifurcation diagrams, stable and unstable branches are blue and red, respectively. In (b) and (c), the branches for synchronous stable steady states are marked by S. These two diagrams are divided into multiple zones based on the number of stable steady states of the whole system (written within circles). BP: Branch Point. LP: Limit Point. The phase diagram of the three-cell system with (d) positive and (e) negative coupling for different $r$ and $p$ values. The number of stable steady states of the system in each parameter zone is written on the diagram.